\documentclass{amsart}
\usepackage[colorlinks=true, pdfstartview=FitV, linkcolor=blue, citecolor=blue, urlcolor=red]{hyperref} 
\usepackage{amssymb,graphicx,color}

\newcommand{\mse}{\ensuremath\operatorname{MSE}}
\newcommand{\expect}{\ensuremath\mathbb{E}}
\newcommand{\vari}{\ensuremath\mathbb{V}}
\newcommand{\bias}{\ensuremath\mathbb{B}}

\newcommand{\colrule}{\hline\hline}
\newcommand{\toprule}{\hline}
\newcommand{\botrule}{\hline}
\newcommand{\tbl}[1]{\caption{#1}}
\newcommand{\subfigure}[1]{\begin{minipage}{.49\textwidth}{#1}\end{minipage}}
\newcommand{\tablefont}{}

\begin{document}

\markboth{P.S.~Chami, B.~Sing and D.~Thomas}{A ratio-product-ratio estimator using auxiliary information}


\title[A ratio-product-ratio estimator using auxiliary information]{A two parameter ratio-product-ratio estimator using auxiliary information}


\author{Peter S. Chami}
\address{The Warren Alpert Medical School of Brown University, 
Box G-A1, 
Providence, Rhode Island 02912, 
USA;
Department of Computer Science, Mathematics and Physics,
Faculty of Pure and Applied Science,
The University of the West Indies, 
P.O. Box 64, Cave Hill,
Bridgetown, St. Michael, BB11000,
BARBADOS W.I. }
\email{peter.chami@cavehill.uwi.edu}
\author{Bernd Sing}
\address{Department of Computer Science, Mathematics and Physics,
Faculty of Pure and Applied Science,
The University of the West Indies, 
P.O. Box 64, Cave Hill,
Bridgetown, St. Michael, BB11000,
BARBADOS W.I.}
\email{bernd.sing@cavehill.uwi.edu}
\urladdr{\htmladdnormallink{http://www.cavehill.uwi.edu/fpas/cmp/staff/bsing/home.html}{http://www.cavehill.uwi.edu/fpas/cmp/staff/bsing/home.html}}
\author{Doneal Thomas}
\address{Department of Computer Science, Mathematics and Physics,
Faculty of Pure and Applied Science,
The University of the West Indies, 
P.O. Box 64, Cave Hill,
Bridgetown, St. Michael, BB11000,
BARBADOS W.I.}
\email{donealt@gmail.com}

\thanks{\textit{Corresponding Author}. Doneal Thomas, \texttt{donealt@gmail.com}}

\date{\today}


\begin{abstract}
We propose a two parameter ratio-product-ratio estimator for a finite population mean in a simple random sample without replacement following the methodology in  Ray et al.~\cite{Ray}, Sahai et al.~\cite{Sahai80,Sahai85} and Singh et al.~\cite{Singh03}.

The bias and mean square error of our proposed estimator are obtained to the first degree of approximation. We derive conditions for the parameters under which the proposed estimator has smaller mean square error than the sample mean, ratio and product estimators.

We carry out an application showing that the proposed estimator outperforms the traditional estimators using groundwater data taken from a geological site in the state of Florida.
\end{abstract}

\keywords{Estimators; Finite population mean; Mean square error; Bias; Relative efficiency}
\subjclass[2010]{62D05; 62G05}

\maketitle

\section{Introduction}

We consider the following setting: For a finite population of size $N$, we are interested in estimating the \emph{population mean} $\bar{Y}$ of the \emph{main variable} $y$ (taking values $y_{i}$ for $i=1,\ldots,N$) from a \emph{simple random sample} of size $n$ (where $n<N$) drawn without replacement. We also know the population mean $\bar{X}$ for the \emph{auxiliary variable} $x$ (taking values $x_{i}$ for $i=1,\ldots,N$). We use the notation $\bar{y}$ and $\bar{x}$ for the \emph{sample means}, which are unbiased estimators of the population means $\bar{Y}$ and $\bar{X}$, respectively. 

We denote the \emph{population variances} of $Y$ and $X$ by
\begin{equation*}
S_{Y}^{2}=\vari{(Y)} =\frac1{N-1}\,\sum\limits_{i=1}^{N}\left( Y_{i}-\bar{Y}\right)^2
\quad\text{and}\quad
S_{X}^{2} =\vari{(X)}=\frac1{N-1}\,\sum\limits_{i=1}^{N}\left( X_{i}-\bar{X}\right)^2,
\end{equation*}
respectively. Furthermore, we define the \emph{coefficient of variation} of $Y$ and $X$ as
\begin{equation*}
C_{Y}=\frac{S_{Y}}{\bar{Y}} \quad\text{and}\quad C_{X}=\frac{S_{X}}{\bar{X}},
\end{equation*}
respectively, and the \emph{coefficient of correlation} $C$ between the two variables as
\begin{equation*}
C=r\cdot\frac{C_{Y}}{C_{X}}
\end{equation*}
where $r=\frac{S_{XY}}{S_{X}S_{Y}}$ denotes \emph{Pearson's correlation coefficient}.

As estimator of the population mean $\bar{Y}$, the usual \emph{sample mean} $\bar{y}$, the \emph{ratio estimator} $\bar{y}_{R}=\left( \bar{y}/\bar{x}\right)\cdot\bar{X}$ and the \emph{product estimator} $\bar{y}_{P}=\left(\bar{x}\cdot\bar{y}\right)/\bar{X}$ are used. Murthy \cite{Murthy} and Sahai at al.~\cite{Sahai80} compared the relative precision of these estimators and showed that the ratio estimator, sample mean and product estimator are most efficient when $C>\frac{1}{2},$ $-\frac{1}{2}\leq C\leq \frac{1}{2}$ and $C<-\frac{1}{2}$, respectively. In other words: when the study variate $y$ and the auxiliary variate $x$ show high positive correlation, then the ratio estimator shows the highest efficiency; when they are highly negative correlated, then the product estimator has the highest efficiency; and when the variables show a weak correlation only, then the sample mean is preferred.

For estimating the population mean $\bar{Y}$ of the main variable, we proposed the following two parameter ratio-product-ratio estimator:
\begin{equation}\label{eq:ratiopro}
\bar{y}_{\alpha,\beta}=\alpha \left[ \frac{\left( 1-\beta \right) \bar{x}+\beta \bar{X}}{\beta \bar{x}+\left( 1-\beta \right) \bar{X}}\right] \bar{y}+\left( 1-\alpha \right) \left[ \frac{\beta \bar{x}+\left( 1-\beta \right) \bar{X}}{\left( 1-\beta \right) \bar{x}+\beta \bar{X}}\right] \bar{y}
\end{equation}
where $\alpha, \beta$ are real constants. Our goal in this article is to derive values for these constants $\alpha, \beta$ such that the bias and/or the mean square error (MSE) of $\bar{y}_{\alpha,\beta}$ is minimal. In fact, in Section~\ref{sec:AOE} we are able to use the two parameters $\alpha$ and $\beta$ to obtain an estimator $\bar{y}^*(C)$ that is (up to first degree of approximation) both unbiased and has minimal MSE.

Note that $\bar{y}_{\alpha,\beta}=\bar{y}_{1-\alpha,1-\beta}$, i.e., the estimator $\bar{y}_{\alpha,\beta}$ is invariant under a point reflection through $(\alpha,\beta)=(\frac12,\frac12)$. In the point of symmetry $(\alpha,\beta)=(\frac12,\frac12)$, the estimator reduces to the sample mean, i.e., we have $\bar{y}_{\frac12,\frac12}=\bar{y}$. In fact, on the whole line $\beta=\frac12$ our proposed estimator reduces to the sample mean estimator, i.e., $\bar{y}_{\alpha,\frac12}=\bar{y}$. Similarly, we get $\bar{y}_{1,0}=\bar{y}_{0,1}=\left( \bar{x}\text{ }\bar{y}\right) /\bar{X}=\bar{y}_{P}$ (product estimator), and $\bar{y}_{0,0}=\bar{y}_{1,1}=\left( \bar{y}\text{ }\bar{X}\right) /\bar{x}=\bar{y}_{R}$ (ratio estimator).

\section{First Degree Approximation to the Bias}\label{sec:bias}

In order to derive the bias of $\bar{y}_{\alpha,\beta}$ up to $O\left( \frac{1}{n}\right)$, we set 
\begin{equation*}
e_{1} =\frac{\bar{y}-\bar{Y}}{\bar{Y}} \quad\text{and}\quad e_{2} =\frac{\bar{x}-\bar{X}}{\bar{X}}.
\end{equation*}%
Thus, we have $\bar{y}=\bar{Y}\left( 1+e_{1}\right)$ and $\bar{x}=\bar{X}\left( 1+e_{2}\right)$, and the relative estimators are given by 
\begin{equation*}
\hat{y} =\frac{\bar{y}}{\bar{Y}}=\left( 1+e_{1}\right)\quad\text{and}\quad\hat{x} =\frac{\bar{x}}{\bar{X}}=\left( 1+e_{2}\right).
\end{equation*}%
Thus, the expectation value of the $e_i$'s is
\begin{equation*}
\expect\left( e_{i}\right) =0\qquad \text{for }\ i=1,2,
\end{equation*}%
and under a simple random sample without replacement, the relative variances are
\begin{equation*}
\vari_{\text{rel}}\left( \bar{y}\right) = \frac{\vari\left(\bar{y}\right)}{\bar{Y}^2} =\expect\left( e_{1}^{2}\right) =\vari\left(e_{1}\right) =\frac{1-f}{n}\left( \frac{S_{Y}}{\bar{Y}}\right) ^{2},
\end{equation*}
and
\begin{equation*}
\vari_{\text{rel}}\left( \bar{x}\right) = \frac{\vari\left(\bar{x}\right)}{\bar{X}^2} = \expect\left( e_{2}^{2}\right) =\vari\left( e_{2}\right) = \frac{1-f}{n}\left( \frac{S_{X}}{\bar{X}}\right)^{2},
\end{equation*}
where $f=\frac{n}{N}$ is the \emph{sampling fraction}.
Also, we have
\begin{equation*}
\expect\left( e_{1}e_{2}\right) =\frac{1-f}{n}\,r\,C_{Y}C_{X},
\end{equation*}
see \cite{Ray,Sahai80,Singh03}.
Furthermore, we note that $\expect\left( e_{1}^{2}e_{2}^{2}\right)= O\left( \frac{1}{n^2}\right)$, and $\expect\left( e_{1}^{i}e_{2}^{j}\right) =0$ when $(i+j)$ is an odd integer.

Now re-expressing \eqref{eq:ratiopro} in terms of $e_i$'s and by substituting $\bar{x}$ and $\bar{y},$ we have
\begin{eqnarray*}
\bar{y}_{\alpha,\beta} &=&\alpha \left[ \frac{1+e_{2}-\beta e_{2}}{ 1+\beta e_{2}}\right] \bar{Y}\left( 1+e_{1}\right) +\left( 1-\alpha \right) \left[ \frac{1+\beta e_{2}}{1+e_{2}-\beta e_{2}}\right] \bar{Y}\left(1+e_{1}\right).
\end{eqnarray*}
In the following, we assume that\footnote{We note that $ \min\left\{\frac1{|\beta|},\frac1{|1-\beta|}\right\}$ attains its maximal value $2$ at $\beta=\frac12$.} $|e_{2}|< \min\left\{\frac1{|\beta|},\frac1{|1-\beta|}\right\}$ so that we can expand $\left( 1+\beta e_{2}\right) ^{-1}$ and $\left( 1+\left( 1-\beta \right)e_{2} \right)^{-1} $ as a series in powers of $e_{2}$ up to $O\left( e_{2}^{3}\right)$. We get
\begin{equation*}
\bar{y}_{\alpha,\beta} = (1+e_1)\,\bar{Y}\cdot\left[1^{}_{}-(1-2\alpha)(1-2\beta)\,e_2+(1-\alpha-\beta)(1-2\beta)\,e^2_2+O(e^3_2)\right].
\end{equation*}

\begin{figure}[tb]
\subfigure{
\resizebox*{5.7cm}{!}{\includegraphics{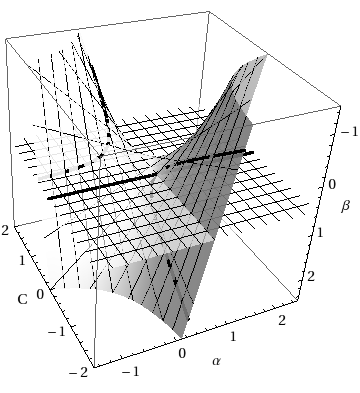}} \begin{center}{\rm (a)}\end{center}}\hfill
\subfigure{
\resizebox*{5.7cm}{!}{\includegraphics{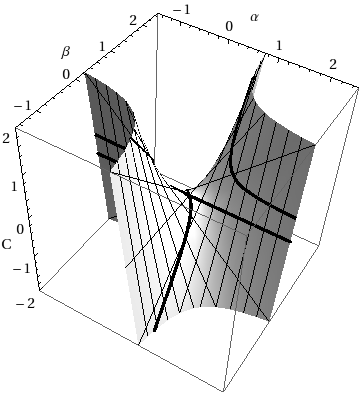}} \begin{center}{\rm (b)}\end{center}}
\caption{\label{fig:biasmse} {\rm (a)}~Surface of ``biasfree estimators'' defined by \eqref{eq:bias} in the parameter space $(\alpha,\beta,C)\in\mathbb{R}^3$.  {\rm (b)}~Surface of ``AOE parameters'' defined by \eqref{eq:sol2}. The points of intersection of the two surfaces (see Section~\ref{sec:AOE}) are drawn as black curves.}
\end{figure}

We assume that the sample is large enough to make $|e_{2}|$ so small that contributions from powers of $e_{2}$ of degree higher than two are negligible, compare \cite{Singh03}. By retaining powers up to $e_{2}^{2}$, we get
\begin{equation}\label{eq:b}
\bar{y}_{\alpha,\beta}-\bar{Y} \approx
\bar{Y}\,\left\{e_1-(1+e_1)\left[(1-2\alpha)(1-2\beta)\,e_2-(1-\alpha-\beta)(1-2\beta)\,e^2_2\right]\right\}
\end{equation}%
Taking expectations on both sides of \eqref{eq:b} and substituting $C=r\frac{C_{Y}}{C_{X}}$, we obtain the bias of $\bar{y}_{\alpha,\beta}$ to order $O\left( n^{-1}\right)$ as
\begin{eqnarray}
\bias\left( \bar{y}_{\alpha,\beta}\right) &=& \expect\left( \bar{y}_{\alpha,\beta}-\bar{Y}\right) \notag\\
&\approx&\frac{1-f}{n}\,(1-2\beta)\,\left[(1-\alpha-\beta)-(1-2\alpha)\,r\,\frac{C_{Y}}{C_X}\right]\,C^2_{X}\,\bar{Y}\notag \\
& = & \frac{1-f}{n}\,(1-2\beta)\,\left[1-\alpha-\beta-(1-2\alpha)\,C\right]\,C^2_{X}\,\bar{Y}\label{eq:bias}
\end{eqnarray}
Equating \eqref{eq:bias} to zero, we obtain
\begin{equation}\label{eq:biasfree}
\beta=\frac12\qquad\text{or}\qquad \beta=1-\alpha-C+2\,\alpha\,C.
\end{equation}
The proposed ratio-product-ratio estimator $\bar{y}_{\alpha,\beta}$, substituted with the values of $\beta$ from \eqref{eq:biasfree}, becomes an (approximately) unbiased estimator for the population mean $\bar{Y}$. In the three-dimensional parameter space $(\alpha,\beta,C)\in\mathbb{R}^3$, these unbiased estimators lie on a plane (in the case $\beta=\frac12$) and on a saddle-shaped surface, see Fig.~\ref{fig:biasmse} {\rm (a)}. Furthermore, as the sample size $n$ approaches the population size $N$, the bias of $\bar{y}_{\alpha,\beta}$ tends to zero, since the factor $(1-f)/n$ clearly tends to zero.

\section{Mean Square Error of \texorpdfstring{${\bar{y}}_{\alpha,\beta}$}{Our Proposed Estimator}}

We calculate the mean square error of $\bar{y}_{\alpha,\beta}$ up to order $O(n^{-1})$ by squaring \eqref{eq:b}, retaining terms up to squares in $e_1$ and $e_2$, and then taking the expectation. This yields 
the first degree approximation of the MSE 
\begin{multline}\label{eq:mse1}
\mse_{1}\left( \bar{y}_{\alpha,\beta}\right) \\ =\frac{1-f}{n}\bar{Y}^{2}\left\{
C_{Y}^{2}+C_{X}^{2}(1-2\alpha)(1-2\beta)\left[(1-2\alpha)(1-2\beta)-2 C\right]\right\}.
\end{multline}
Taking the gradient $\nabla=\left(\frac{\partial}{\partial\alpha},\frac{\partial}{\partial\beta}\right)$ of \eqref{eq:mse1}, we get
\begin{equation}\label{eq:grad}
 \nabla\mse_{1}\!\left( \bar{y}_{\alpha,\beta}\right) =4\,\frac{1-f}{n}\bar{Y}^2\,C_{X}^{2}\,\left[(1-2\alpha)\,(1-2\beta)-C\right]\ \left(1-2\beta^{}_{},1-2\,\alpha^{}_{}\right).
\end{equation}
Setting \eqref{eq:grad} to zero to obtain the critical points, we obtain the following solutions: 
\begin{equation}\label{eq:sol1}
\alpha =\frac{1}{2}\quad\text{and}\quad\beta =\frac{1}{2},
\end{equation}
or
\begin{equation}\label{eq:sol2}
C=(1-2\alpha)\,(1-2\beta).
\end{equation}
One can check that the critical point in \eqref{eq:sol1} is a saddle point unless $C=0$, in which case we get a local minimum. However, the critical points determined by \eqref{eq:sol2} are always local minima; for a given $C$, equation \eqref{eq:sol2} is the equation of a hyperbola symmetric through $(\alpha,\beta)=(\frac12,\frac12)$. Thus, in the three-dimensional parameter space $(\alpha,\beta,C)\in\mathbb{R}^3$, the estimators with minimal MSE (or better, minimal first approximation to the MSE, see calculation in \eqref{eq:minimalMSE} below) lie on a saddle-shaped surface, see Fig.~\ref{fig:biasmse} {\rm (b)}.

We now calculate the minimal value of the MSE. Substituting \eqref{eq:sol1} into the estimator $\bar{y}_{\alpha,\beta}$, yields the unbiased estimator $\bar{y}$ (sample mean) of the population mean $\bar{Y}$. Thus, we arrive at the mean square error of the sample mean:
\begin{equation*}
\mse\left(\bar{y}_{\frac{1}{2},\frac{1}{2}}\right) = \mse\left( \bar{y}\right) =\frac{1-f}{n}%
\bar{Y}^{2}C_{Y}^{2} = \frac{1-f}{n}S_{Y}^{2}
\end{equation*}
By substituting \eqref{eq:sol2} into the estimator, an \emph{asymptotically optimum estimator} (AOE) $\bar{y}_{\alpha,\beta}^{\left(o\right) }$ is found. For the first degree approximation of the MSE, we find
\begin{equation}\label{eq:minimalMSE}
\mse_1\!\left( \bar{y}_{\alpha,\beta}^{\left( o\right) }\right) =\frac{1-f}{n}\bar{Y}^{2}\left( C_Y^2-C^2_{}\,C_X^2\right)=\frac{1-f}{n}S_{Y}^{2}\left( 1-r^{2}\right),
\end{equation}
i.e., the same minimal mean square error as found in \cite{Ray,Sahai80,Sahai85,Singh03}.

\section{Comparison of MSE's and Choice of Parameters}

Here we compare $\mse_{1}\left( \bar{y}_{\alpha,\beta}\right)$ in \eqref{eq:mse1} with the MSE of the product, ratio and sample mean estimators respectively. It is known \cite{Ray,Sahai80} that
\begin{equation}\label{eq:sample}
\mse\left( \bar{y}\right) =\vari\left( \bar{y}\right) =\frac{1-f}{n}\bar{Y}^{2}C_{Y}^{2},
\end{equation}
\begin{equation}\label{eq:ratio}
\mse_1\left( \bar{y}_{R}\right) =\frac{1-f}{n}\bar{Y}^{2}\{C_{Y}^{2}+C_{X}^{2}\left( 1-2C\right) \},
\end{equation}
and
\begin{equation}\label{eq:product}
\mse_1\left( \bar{y}_{P}\right) =\frac{1-f}{n}\bar{Y}^{2}\{C_{Y}^{2}+C_{X}^{2}\left( 1+2C\right) \}.
\end{equation}

\subsection{Comparing the MSE of the product estimator to our proposed estimator}

From \cite{Ray,Sahai80,Sahai85,Singh03}, we know that for $C<-\frac{1}{2}$, the product estimator is preferred to the sample mean and ratio estimators. Therefore, we seek a range of $\alpha$ and $\beta$ values where our proposed estimator $\bar{y}_{\alpha,\beta}$ has smaller MSE than the product estimator.

From \eqref{eq:product} and \eqref{eq:mse1}, the following expression can be verified:
\begin{equation*}
\mse_1\left( \bar{y}_{P}\right) -\mse_{1}\left( \bar{y}_{\alpha,\beta}\right) =4\,\frac{1-f}{n}\bar{Y}^{2}\, C_{X}^{2} \left[1+2\alpha\beta-\alpha-\beta\right]\,\left[C-(2\alpha\beta-\alpha-\beta)\right]
\end{equation*}
which is positive if
\begin{equation}\label{eq:diff1}
\left[1+2\alpha\beta-\alpha-\beta\right]\,\left[C-(2\alpha\beta-\alpha-\beta)\right]>0.
\end{equation}%
We obtain the following two cases
\begin{enumerate}
\item\label{it:prod1} $C>2\alpha\beta-\alpha-\beta>-1$ (if both factors in \eqref{eq:diff1} are positive), or
\item\label{it:prod2} $C<2\alpha\beta-\alpha-\beta<-1$ (if both factors in \eqref{eq:diff1} are negative).
\end{enumerate}

Noting that we are only interested in the case $C<-\frac{1}{2}$, we get from \ref{it:prod1}
\begin{equation*}
-\frac12>C>2\alpha\beta-\alpha-\beta>-1.
\end{equation*}%
We note that this implies $-1<C<-\frac12$, and the range for $\alpha $ and $\beta $ where these inequalities hold are explicitly given by the following two cases:
\begin{itemize}
\item If $\beta <\frac{1}{2}$, then $\frac{\beta+C}{2\beta -1}<\alpha <\frac{\beta -1}{2\beta -1}$, and
\item if $\beta >\frac{1}{2}$, then $\frac{\beta-1}{2\beta -1}<\alpha <\frac{\beta+C}{2\beta -1}$.
\end{itemize}
For any given $C$, we again note that the two regions determined here are symmetric through $(\alpha,\beta)=(\frac12,\frac12)$. We also note that the parameters $(\alpha,\beta)$ which give an AOE (see \eqref{eq:sol2}), which for a fixed $C$ lie on a hyperbola, are contained in these region.   

In the case \ref{it:prod2}, where $C<-1$ (and therefore automatically $C<-\frac12$), the following range for $\alpha$ and $\beta$ can be found:
\begin{itemize}
\item If $\beta <\frac{1}{2}$, then $\frac{\beta-1}{2\beta -1}<\alpha <\frac{\beta+C}{2\beta-1}$, and
\item if $\beta >\frac{1}{2}$, then $\frac{\beta+C}{2\beta -1}<\alpha <\frac{\beta-1}{2\beta-1}$.
\end{itemize}
The same remark as in the previous case apply. Furthermore, note that for $C=-1$, the product estimator attains the same minimal MSE as our proposed estimator $\bar{y}_{\alpha,\beta}$ on the hyperbola given by \eqref{eq:minimalMSE}. In Fig.~\ref{fig:msedata} {\rm (a)} we show the region in parameter space $(\alpha,\beta,C)\in\mathbb{R}^3$ calculated here and in the next two sections where the proposed estimator works better than the three traditional estimators.

\subsection{Comparing the MSE of the ratio estimator to our proposed estimator}\label{sec:ratio}

For $C>\frac{1}{2},$ the ratio estimator is used instead of the sample mean or product estimator \cite{Ray,Sahai80,Sahai85,Singh03}. As a result, we are concern with a range of plausible values for $\alpha $ and $\beta $, where $\bar{y}_{\alpha,\beta}$ works better than the ratio estimator.

Taking the difference of \eqref{eq:ratio} and \eqref{eq:mse1}, we have
\begin{equation*}
\mse_1\left( \bar{y}_{R}\right) -\mse_{1}\left( \bar{y}_{\alpha,\beta}\right) =4\,\frac{1-f}{n}\bar{Y}^{2}\, C_{X}^{2} \left[2\alpha\beta-\alpha-\beta\right]\,\left[C-1-(2\alpha\beta-\alpha-\beta)\right]
\end{equation*}
which is positive if
\begin{equation*}
\left[2\alpha\beta-\alpha-\beta\right]\,\left[C-1-(2\alpha\beta-\alpha-\beta)\right]>0.
\end{equation*}
Therefore, 
\begin{enumerate}
\item\label{it:ratio1} $C-1>2\alpha\beta-\alpha-\beta>0$, or
\item\label{it:ratio2} $C-1<2\alpha\beta-\alpha-\beta<0$. 
\end{enumerate}
Hence, from solution \ref{it:ratio1}, where $C>1$, we have
\begin{itemize}
\item If $\beta <\frac{1}{2}$, then $\frac{\beta+C-1}{2\beta -1}<\alpha <\frac{\beta}{2\beta -1}$, and
\item if $\beta >\frac{1}{2}$, then $\frac{\beta}{2\beta -1}<\alpha <\frac{\beta+C-1}{2\beta -1}$.
\end{itemize}
Also, from solution \ref{it:ratio2}, where $\frac{1}{2}<C<1$, we obtain
\begin{itemize}
\item If $\beta <\frac{1}{2}$, then $\frac{\beta}{2\beta -1}<\alpha <\frac{\beta+C-1}{2\beta -1}$, and
\item if $\beta >\frac{1}{2}$, then $\frac{\beta+C -1}{2\beta -1}<\alpha <\frac{\beta}{2\beta -1}$.
\end{itemize}

\subsection{Comparing the MSE of the sample mean to our proposed estimator}

Finally, we compare the $\mse\left( \bar{y}\right) $ to our proposed estimator, $\mse\left( \bar{y}_{\alpha,\beta}\right)$. From \cite{Ray,Sahai80,Sahai85,Singh03}, we know that sample mean estimator is preferred for $-\frac12\le C\le\frac{1}{2}$.

Taking the difference of \eqref{eq:sample} and \eqref{eq:mse1}, we get:
\begin{equation*}
\mse\left( \bar{y}\right) -\mse_1\left( \bar{y}_{\alpha,\beta}\right) =\frac{1-f}{n}\bar{Y}^{2}C_{X}^{2}\,(1-2\alpha)(1-2\beta)\left\{2C-(1-2\alpha)(1-2\beta)\right\}
\end{equation*}
which is positive if 
\begin{equation*}
(1-2\alpha)(1-2\beta)\left\{2C-(1-2\alpha)(1-2\beta)\right\} >0.
\end{equation*}
Therefore, either
\begin{enumerate}
\item\label{it:sample1} $\alpha >\frac{1}{2}$, $\beta >\frac{1}{2}$ and $C>\frac12 (1-2\alpha)(1-2\beta)$, or
\item\label{it:sample2} $\alpha <\frac{1}{2}$, $\beta >\frac{1}{2}$ and $C<\frac12 (1-2\alpha)(1-2\beta)$, or 
\item\label{it:sample3} $\alpha >\frac{1}{2}$, $\beta <\frac{1}{2}$ and $C<\frac12 (1-2\alpha)(1-2\beta)$, or 
\item\label{it:sample4} $\alpha <\frac{1}{2}$, $\beta <\frac{1}{2}$ and $C>\frac12 (1-2\alpha)(1-2\beta)$.
\end{enumerate}
Combining these with the condition $-\frac12\le C\le\frac{1}{2}$, we get the following explicit ranges:
\begin{itemize}
\item If $0<C\le\frac12$ and $\beta >\frac{1}{2}$, then $ \frac12<\alpha <\frac{2\beta+2C-1}{2(2\beta-1)}$  (from \ref{it:sample1}), and
\item if $0<C\le\frac12$ and $\beta <\frac{1}{2}$, then $\frac{2\beta+2C-1}{2\left( 2\beta -1\right) }<\alpha <\frac{1}{2}$ (from \ref{it:sample4}), and
\item if $-\frac12\le C<0$ and $\beta >\frac{1}{2}$, then $\frac{2\beta+2C-1}{2\left( 2\beta -1\right) }<\alpha <\frac{1}{2}$ (from \ref{it:sample2}), and
\item if $-\frac12\le C<0$ and $\beta <\frac{1}{2}$, then $ \frac{1}{2}<\alpha <\frac{2\beta+2C-1}{2\left( 2\beta -1\right) }$ (from \ref{it:sample3}).
\end{itemize}
We note that the case $C=0$ implies $r=0$, and thus the sample mean estimator is the estimator with minimal MSE (and, as already noted, $\bar{y}=\bar{y}_{\frac12,\frac12}$).

\begin{figure}[tb]
\subfigure{
\resizebox*{5.7cm}{!}{\includegraphics{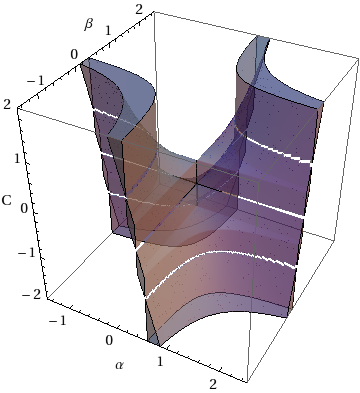}} \begin{center}{\rm (a)}\end{center} }\hfill
\subfigure{
\resizebox*{5.7cm}{!}{\includegraphics{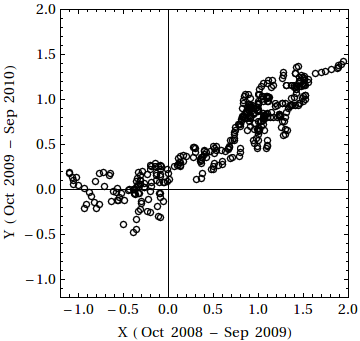}}\linebreak \begin{center}{\rm (b)}\end{center} }
\caption{\label{fig:msedata} {\rm (a)}~Region of the parameter space $(\alpha,\beta,C)\in\mathbb{R}^3$ where our proposed estimator $\bar{y}_{\alpha,\beta}$ has lower MSE than the traditional estimators. {\rm (b)}~Scatterplot of study and auxiliary variables for the groundwater data studied in Section~\ref{sec:app}.}
\end{figure}

In Fig.~\ref{fig:msedata} {\rm (a)} we show the region in parameter space $(\alpha,\beta,C)\in\mathbb{R}^3$  where the proposed estimator $\bar{y}_{\alpha,\beta}$ works better than the three traditional estimators. Note that the surface of ``AOE parameters'' in Fig.~\ref{fig:biasmse} {\rm (b)} is a subset of this region, except for the values $C=0$, $C=-1$ and $C=+ 1$ for which our proposed estimator only works as well as\footnote{We also remark that the points $(1,1,1)$ and $(0,0,1)$ (note that $\bar{y}_{1,1}=\bar{y}_{0,0}=\bar{y}_{R}$), $(0,1,-1)$ and $(1,0,-1)$ (note that $\bar{y}_{1,0}=\bar{y}_{0,1}=\bar{y}_{R}$) as well as the line $(\alpha,\frac12,0)$ (note that $\bar{y}_{\alpha,\frac12}=\bar{y}$) belong to the surface of ``AOE parameters'' in Fig.~\ref{fig:biasmse} {\rm (b)}.} the sample mean, product and ratio estimator, respectively.

\section{Unbiased AOE}\label{sec:AOE}

Combining \eqref{eq:biasfree} and \eqref{eq:sol2}, we can calculate the parameters $\alpha$ and $\beta$ where our proposed estimator becomes -- at least up to first approximation -- an unbiased AOE. We obtain a line with
\begin{equation*}
\beta=\frac12 \quad\text{and}\quad C=0,
\end{equation*}
(recall that on this line our estimator always reduces to the sample mean estimator) or a ``curve'' $(\alpha^*(C),\beta^*(C),C)\in\mathbb{R}$ in the parameter space with
\begin{equation}\label{eq:opt}
\alpha^*(C)=\frac12\left(1\pm\sqrt{\frac{C}{2C-1}}\right) \quad\text{and}\quad \beta^*(C)=\frac12\left(1\pm\sqrt{C(2C-1)}\right).
\end{equation}
We note that the parametric ``curve'' in \eqref{eq:opt} is only defined for $C\le 0$ or $C>\frac12$ -- in fact, this parametric ``curve'' is three hyperbolas.  The surface of ``biasfree estimator parameters'' in Fig.~\ref{fig:biasmse} {\rm (a)} and the surface of ``AOE parameters'' in Fig.~\ref{fig:biasmse} {\rm (b)} only intersect in these three hyperbolas and the line $\beta=\frac12$ and $C=0$. In the region $0< C\le\frac12$ of the parameter space $(\alpha,\beta,C)\in\mathbb{R}^3$, we have the common situation where minimising MSE comes with a trade-off in bias. The curves of intersection are included in Fig.~\ref{fig:biasmse}. Explicitly, our proposed estimator using the values in \eqref{eq:opt} is given by
\begin{multline}\label{eq:unbiasedAOE}
\bar{y}^*(C)=\bar{y}_{\alpha^*(C),\beta^*(C)} =\\
\frac{2\,(C+1)\,\bar{X}^2-2\,(C-1)\,\bar{x}^2+\left(2\,C^2-C-1\right)\,\left(\bar{X}-\bar{x}\right)^2}{4\,\bar{X}\,\bar{x}-\left(2\,C^2-C-1\right)\,\left(\bar{X}-\bar{x}\right)^2}\,\bar{y}.
\end{multline} 
At first it might seem surprising that this estimator $\bar{y}^*(C)$ is also defined in\footnote{The denominator vanishes if $C=\frac14\left(1\pm\frac{\sqrt{9\,(\bar{X}-\bar{x})^2+32\,\bar{X}\,\bar{x}}}{\bar{X}-\bar{x}}\right)$.} the region $0< C\le\frac12$. However, one can also let the parameters $(\alpha,\beta)$ in the definition of our proposed estimator $\bar{y}_{\alpha,\beta}$ in \eqref{eq:ratiopro} be complex numbers -- but such that we still get a real estimator.  One can check that $\alpha^*(C)$ and $\beta^*(C)$ in \eqref{eq:opt} for $0<C<\frac12$ have this property. 

Furthermore, we can check that the first degree of approximation of the bias and MSE of $\bar{y}^*(C)$ are given by
\begin{equation*}
\bias_1\left(\bar{y}^*(C)\right)=0 \qquad\text{and}\qquad \mse_1\left(\bar{y}^*(C)\right) = \frac{1-f}{n}S_{Y}^{2}\left( 1-r^{2}\right)
\end{equation*}
(compare \eqref{eq:bias} and \eqref{eq:minimalMSE}). Thus, the estimator $\bar{y}^*(C)$ of \eqref{eq:unbiasedAOE} is an unbiased AOE.

One might also ask whether inside $0< C\le\frac12$ there is a choice of real parameters $(\alpha,\beta)\in\mathbb{R}^2$ such that we get an AOE with small bias. Using \eqref{eq:sol2} in \eqref{eq:bias}, we get the first degree approximation of the bias of an AOE
\begin{equation*}
\bias_1\left(\bar{y}^{(o)}_{\alpha,\beta}\right) = \frac{1-f}{n}C_X^2 \bar{Y}\,\frac12\left[C\left(1-2\,C\right)+(1-2\beta)^2\right].
\end{equation*}
From this expression (and the constraint \eqref{eq:sol2}) it is clear that the bias can only be made zero if $C\le0$ or $C\ge\frac12$. Otherwise, there is always a positive contribution coming from the term $C(1-2C)$ that does not vanish no matter what we choose for $\beta\in\mathbb{R}$. In fact, it looks as if the choice $\beta=\frac12$ always yields the least possible bias, however two remarks are in order here: Firstly, given \eqref{eq:sol2} and unless $C=0$, we can only let $\beta$ be close to $\frac12$ and choose $\alpha$ accordingly (the absolute value $|\alpha|$ is then large). Secondly, we already noted that $\bar{y}_{\alpha,\frac12}=\bar{y}$, and the MSE for the sample mean estimator is $\frac{1-f}{n} S^2_Y$, not $\frac{1-f}{n} S^2_Y(1-r^2)$ as for an AOE. We have arrived here at a point where the first degree approximation to bias and MSE breaks down. To find a choice of real parameters for given $C$ with minimal MSE and least bias, higher degrees of approximation would have to be considered.

\section{Application and Conclusion}\label{sec:app}

Using data taken from the Department of the Interior, United States Geological Survey (USGS) \cite{USGS}, site number 02290829501 (located in Florida), a comparison of our proposed estimator $\bar{y}_{\alpha,\beta}$ to the traditional estimators was carried out. The study variables (denoted by $Y$) are taken to be the maximum daily values (in feet) of groundwater at the site for the period October 2009 to September 2010. The auxiliary variables (denoted $X$) are taken as the maximum daily values (in feet) of groundwater for the period October 2008 to September 2009. Our goal is to estimate the true average maximum daily groundwater $\bar{Y}$ for the period October 2009 to September 2010.

The questions we ask are: How many units of groundwater must be taken from the population $Y$ to estimate the population mean $\bar{Y}$ within $d=10\%$ at a $90\%$ confidence level $(\alpha=0.10)$? And how well do the estimators perform given this data set with auxiliary information for the calculated sample size $n$?

Using the entire data set, we calculate the following statistics: 
$\bar{Y}\approx0.5832$, $\bar{X}\approx0.6277$, $S_{Y}\approx0.4480$, $S_{X}\approx0.7222$, $r\approx0.9125$, $C_{Y}\approx0.7681$, $C_{X}\approx1.1504$ and $C\approx0.6092$.
A scatterplot of the data set is shown in Fig.~\ref{fig:msedata} {\rm (b)}, which adds emphasis to the positive measure of association between the study variable $Y$ and the auxiliary variable $X$.

One should note that the value of $C\approx0.6092$ lies in the interval $(\frac{1}{2},1)$, so we choose values of $\alpha$ and $\beta$ from Section~\ref{sec:ratio}. Indeed, we use \eqref{eq:opt} and choose $\beta=\beta^*(0.6092)\approx 0.3176$ and $\alpha=\alpha^*(0.6092)\approx -0.3349$. Note that $\beta=0.3176$ yields $(-0.8704)<\alpha<0.200573$  in Section~\ref{sec:ratio}. Using the notation of Section~\ref{sec:AOE}, we also note that $\bar{y}_{-0.3349,0.3176} = \bar{y}^*(0.6092)$.

\subsection{Calculating the sample size \texorpdfstring{$n$}{n}}

To estimate the population mean amount of groundwater recorded for the state of Florida from October 2009 to September 2010, a sample of size $n$ is drawn from the population of size $N=365$ according to the simple random sampling without replacement \cite{Cochran}. A first approximation to this sample size needed is the (infinite population) value   
\begin{equation*}
n_{0}=\frac{Z_{\frac{\alpha}{2}}^{2}\sigma^{2}}{d^{2}},
\end{equation*}
where $d$ is the chosen margin of error from the estimate of $\bar{Y}$, and $Z_{\frac{\alpha}{2}}$ is a standard normal variable with tail probability of $\frac{\alpha}{2}$. Accounting for the finite population size $N$, we obtain the sample size
\begin{equation*}
n=\frac1{\frac1{n_0}+\frac1{N}}.
\end{equation*}
In general, the true value of $\sigma^{2}$ is unknown but can be estimated using its consistent estimator $s^{2}$. However, in our case $\sigma^{2}$ is calculated from the population and is given as $S_{Y}^{2}\approx0.2006$. Therefore, with $\alpha=0.10$, $d=10\%\text{ of }\bar{Y}$ (i.e., $d\approx 0.0583$) and $Z_{0.05}\approx1.6449$, the sample size can be calculated as  $n_{0}\approx\frac{(1.6449)^{2}\cdot 0.2006}{(0.0583)^{2}}\approx 159.59$, rounding up gives $n_{0}=160$. Therefore, we get $n\approx\frac1{\frac1{160}+\frac1{365}}\approx 111.23$, and thus take $n=112$.

\subsection{Relative efficiencies}

Table~\ref{tab:eff} shows the relative efficiencies of the traditional estimators (sample mean $\bar{y}$, ratio $\bar{y}_{R}$ and product $\bar{y}_{P}$ estimators) and our proposed two parameter ratio-product-ratio estimator $\bar{y}_{\alpha,\beta}$ for $(\alpha,\beta)=(-0.3349,0.3176)$. We note that with this choice of parameters, the estimator is an (unbiased) AOE, namely $\bar{y}_{-0.3349,0.3176}=\bar{y}^*(0.6092)$. The table shows that our two parameter ratio-product-ratio estimator dominates the traditional estimators in the sense that it has the highest efficiency. 

We can also observe that in the computation of the relative efficiency, the specification of the sample size $n$ is not important since the finite population correction factor $(\frac{1-f}{n})$ cancels out (however, this would not be the case for higher degrees of approximation). 

\begin{table}
\tbl{Relative efficiencies comparisons.}
{\begin{tabular}{cccc}
\toprule
$\frac{\mse(\bar{y})}{\mse(\bar{y})}$ & $\frac{\mse(\bar{y})}{\mse_1(\bar{y}_{R})}$ & $\frac{\mse(\bar{y})}{\mse_1(\bar{y}_{P})}$ & $\frac{\mse(\bar{y})}{\mse_1\left(\bar{y}^*(0.6092)\right)}$\rule[-2.5mm]{0mm}{0mm}\rule[4mm]{0mm}{0mm}\\
\colrule
$100\%$ & $196.11\%$ & $16.73\%$ & $597.28\%$\\
\botrule
\end{tabular}}
\label{tab:eff}
\end{table}

\subsection{Constructing a \texorpdfstring{$90\%$}{90\%} confidence interval for \texorpdfstring{$\bar{Y}$}{Y} using \texorpdfstring{${\bar{y}}_{\alpha,\beta}$}{our proposed estimator}}

Constructing a $90\%$ confidence interval, the following formulation can be used (similar formulae hold for all estimators discussed here), see \cite{Cochran}:
\begin{equation*}
\left(\bar{y}^*(0.6092)\pm Z_{0.05}\sqrt{\frac{S^2_Y}{n}}\cdot\sqrt{\frac{N-n}{N-1}}\right).
\end{equation*}
The factor $\sqrt{(N-n)/(N-1)}$ is the \emph{finite population correction}.

Of course, by the choice of the sample size $n=112$, we get a margin of error of approximately $0.1\cdot \bar{Y}\approx 0.0583$ here; more precisely, the calculation using the above formula yields $\left(\bar{y}^*(0.6092)\pm 1.6449\cdot 0.0423\cdot 0.8337\right)=\left(\bar{y}^*(0.6092) \pm 0.0580\right)$.

\subsection{Comparison of estimators}

To compare the proposed estimator with the traditional ones, we selected $10\,000$ times a sample of size $n=112$ and calculate the estimators from it. We note that  there are $\binom{365}{112}\approx 2.5\cdot10^{96}$ possibilities to choose $112$ data points out of a total $365$ without replacement.

In Table~\ref{tab:bias} we show the relative position of the estimators with respect to the population mean $\bar{Y}$. In the $10\,000$ simulations, our proposed estimator outperformed the traditional estimators on $4\,701$ occasions. The ratio estimator, the suggested estimator for this value of $C$ by~\cite{Murthy}, performs better than our proposed estimator $2\,978$ times (in these cases it is actually the best of the studied estimators; note that the ratio estimator is the worst $1\,349$ times).  

In Table~\ref{tab:est}, we compare the estimators by looking at the following criteria: The \emph{coverage probability} is the proportion of the $90\%$ confidence interval covering the population mean $\bar{Y}$; as expected, the usual mean sample estimator yields around $90\%$, while the ratio estimator and our proposed estimator yield much higher values -- in this simulation, all intervals calculated from our proposed estimator cover $\bar{Y}$. For those $90\%$ confidence intervals that do not cover $\bar{Y}$, we check whether they lie to the left (\emph{negative bias}) or to the right (\emph{positive bias}) of $\bar{Y}$. We also state the statistical information \emph{lower} and \emph{upper quartile} and \emph{median}, we get from the $10\,000$ simulations; directly below this we show violin plots for the estimators (the dashed line indicates the value $\bar{Y}$, the dotted lines indicate the $90\%$ confidence interval). In the violin plot, we see that the values obtained by our proposed estimators yield a narrow normal distribution around the true value (skewness is $0.0046$, kurtosis is $2.9926$), while the product estimator gives a spread-out distribution and the (traditionally preferred) ratio estimator gives a skewed distribution (skewness is $0.5230$, kurtosis is $3.4253$). Finally, we compare the values of the MSEs; the experimental values obtained agree with the theoretical values listed in Table~\ref{tab:eff}.       

We infer that our proposed estimator is more efficient and robust than the traditional sample mean, ratio and product estimators.

\begin{table}
\tbl{Comparison of the estimators according to the absolute deviation from the population mean $\bar{Y}$ (in $10\,000$ simulations).}
{\begin{tabular}{|c@{ $<$ }c@{ $<$ }c@{ $<$ }c|r|} \toprule
\multicolumn{4}{|c|}{position} & counts \\ \colrule
$\left|\bar{y}^*(0.6092)-\bar{Y}\right|$ & $\left|\bar{y}^{}_{R}-\bar{Y}\right|$ & $\left|\bar{y}^{}_{}-\bar{Y}\right|$ & $\left|\bar{y}^{}_{P}-\bar{Y}\right|$ & 3\,173\\[.5mm]
$\left|\bar{y}^{}_{R}-\bar{Y}\right|$ & $\left|\bar{y}^*(0.6092)-\bar{Y}\right|$ & $\left|\bar{y}^{}_{}-\bar{Y}\right|$ & $\left|\bar{y}^{}_{P}-\bar{Y}\right|$ & 2\,978\\[.5mm]
$\left|\bar{y}^*(0.6092)-\bar{Y}\right|$ & $\left|\bar{y}^{}_{}-\bar{Y}\right|$ & $\left|\bar{y}^{}_{R}-\bar{Y}\right|$ & $\left|\bar{y}^{}_{P}-\bar{Y}\right|$ & 1\,528\\[.5mm]
$\left|\bar{y}^{}_{}-\bar{Y}\right|$ & $\left|\bar{y}^*(0.6092)-\bar{Y}\right|$ & $\left|\bar{y}^{}_{R}-\bar{Y}\right|$ & $\left|\bar{y}^{}_{P}-\bar{Y}\right|$ & 972\\[.5mm]
$\left|\bar{y}^{}_{P}-\bar{Y}\right|$ & $\left|\bar{y}^{}_{}-\bar{Y}\right|$ & $\left|\bar{y}^*(0.6092)-\bar{Y}\right|$ & $\left|\bar{y}^{}_{R}-\bar{Y}\right|$ & 766\\[.5mm]
$\left|\bar{y}^{}_{}-\bar{Y}\right|$ & $\left|\bar{y}^*(0.6092)-\bar{Y}\right|$ & $\left|\bar{y}^{}_{P}-\bar{Y}\right|$ & $\left|\bar{y}^{}_{R}-\bar{Y}\right|$ & 308\\[.5mm]
$\left|\bar{y}^{}_{}-\bar{Y}\right|$ & $\left|\bar{y}^{}_{P}-\bar{Y}\right|$ & $\left|\bar{y}^*(0.6092)-\bar{Y}\right|$ & $\left|\bar{y}^{}_{R}-\bar{Y}\right|$ & 275\\
\botrule
\end{tabular}}
\label{tab:bias}
\end{table}

\begin{table}
\tbl{Comparison of the estimators in $10\,000$ simulations. See text for details.}
{\begin{tabular}{|c|c|cc|ccc|} \toprule
estimator & coverage & neg.~bias & pos.~bias & lo.~quart. & median & up.~quart. \\ \colrule 
$\bar{y}^{}_{}$                    & $89.86\%$ & $5.14\%$ & $5.00\%$ & $0.5583$ & $0.5820$ & $0.6062$  \\ 
$\bar{y}^{}_{R}$                 & $97.01\%$ &  $0.12\%$ & $2.87\%$ & $0.5670$ & $0.5831$ & $0.6017$ \\ 
$\bar{y}^{}_{P}$                 & $49.59\%$ & $25.44\%$ & $24.97\%$ & $0.5240$ & $0.5809$ & $0.6410$  \\ 
$\bar{y}^{*}(0.6092)$         & $100.00\%$ & $0.00\%$ & $0.00\%$ & $0.5732$ & $0.5829$ & $0.5928$\rule[-1.8mm]{0mm}{0mm} \\ \botrule
\end{tabular}}
\begin{minipage}[c]{\textwidth}
\begin{minipage}[c]{.01\textwidth}
\mbox{ }
\end{minipage}
\begin{minipage}[c]{.4\textwidth}\begin{flushleft}
\tablefont
{\begin{tabular}{|c|cc|} \toprule
estimator & MSE & $\frac{\mse(\bar{y}^{}_{})}{\mse(\text{est})}$\rule[-2.2mm]{0mm}{0mm}\rule[3.8mm]{0mm}{0mm} \\ \colrule 
$\bar{y}^{}_{}$                    & $0.00127$ & $100\%$ \\ 
$\bar{y}^{}_{R}$                 & $0.00068$ & $186.68\%$\\ 
$\bar{y}^{}_{P}$                 & $0.00755$ & $16.79\%$ \\
$\bar{y}^{*}(0.6092)$         & $0.00021$ & $598.93\%$\rule[-1.8mm]{0mm}{0mm}\\ \botrule
\end{tabular}}
\end{flushleft}
\end{minipage}\hfill
\begin{minipage}[c]{.52\textwidth}
\includegraphics[scale=0.45]{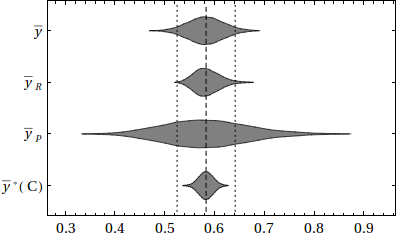}
\end{minipage}
\end{minipage}
\label{tab:est}
\end{table}

\subsection{Outlook} 

Several authors have proposed efficient estimators using auxiliary information. For example,  Srivastava and Reddy \cite{Srivastava67,Reddy} consider a generalisation to the product and ratio estimator given by $\bar{y}^{(k)} = \bar{y} \left({\bar{x}}/{\bar{X}}\right)^k$; Reddy \cite{Reddy} also introduces the estimator $\bar{y}_k={\bar{y}\,\bar{X}}/\left({\bar{X}+k\,(\bar{x}-\bar{X})}\right)$; in Sahai et al.~\cite{Sahai80} the estimator $\bar{y}_{k\,t}=\bar{y}\,\left(2-(\bar{x}/\bar{X})^k\right)$ (where ``t'' stands for ``transformed'') is considered; and Singh et al.~\cite{Singh03} introduce a certain class of ``ratio-product'' estimators given by  $\bar{y}_{RP}(k)=\bar{y} \left(k\cdot (\bar{X}/\bar{x})+(1-k)\cdot(\bar{x}/\bar{X})\right)$. Choosing appropriate parameters $k$ for these estimators and calculating the first degree approximation of the MSE, one can show that  
\begin{multline*}
\mse_1\left(\bar{y}^{(-C)} \right) = \mse_1\left(\bar{y}_C\right) = \mse_1\left(\bar{y}_{C\,t}\right) = \mse_1\left(\bar{y}_{RP}({\scriptstyle\frac{C+1}2}) \right) \\ = \frac{1-f}{n}S_{Y}^{2}\left( 1-r^{2}\right).
\end{multline*}
Thus, these estimators and our proposed estimator (see \eqref{eq:minimalMSE}) are equally efficient up to the first degree of approximation. Indeed, all these estimators give similar results as our proposed estimator in the above application, see Table~\ref{tab:est2}. Comparing the first degree of approximation of the bias (doing calculations as in Section~\ref{sec:bias}) reveals why our unbiased AOE $\bar{y}^*(C)$ and Reddy's $\bar{y}_C$ behave similar -- they are both unbiased AOEs:
\begin{multline*}
\bias(\bar{y}) = \bias_1\left(\bar{y}^*(C)\right) = \bias_1\left(\bar{y}_C\right) = 0, \quad \bias_1\left(\bar{y}_R\right)=\frac{1-f}{n} (1-C) C^2_X \bar{y}, \\ \bias_1\left(\bar{y}_P\right)=\frac{1-f}{n} C C^2_X \bar{y}, \quad \bias_1\left(\bar{y}^{(-C)}\right)=\frac{1-f}{n} \frac{C (1-C)}2 C^2_X \bar{y},\\ \bias_1\left(\bar{y}_{C\,t}\right)=\frac{1-f}{n} \frac{C (1-3\,C)}2 C^2_X \bar{y},  \  \bias_1\left(\bar{y}_{RP}({\scriptstyle\frac{C+1}2})\right)=\frac{1-f}{n} \frac{(1+2C) (1-C)}2 C^2_X \bar{y}.
\end{multline*}
(With $C=0.6092$, only $\bar{y}_{C\,t}$ is negatively biased, compare the quartiles and the box plot in Table~\ref{tab:est2}). For our proposed estimator $\bar{y}_{\alpha,\beta}$, we are able to use the two parameters $\alpha$ and $\beta$ to obtain an estimator $\bar{y}^*(C)$ that is up to first degree of approximation both unbiased and has minimal MSE. A thorough comparison of these estimators involving higher degrees of approximation of MSE and bias as well as accompanying simulations might be desirable, e.g., to find the estimator that behaves well if the coefficient of correlation $C$ is not known exactly. 

\begin{table}
\tbl{Comparison of AOEs in $10\,000$ simulations. See text for details.}
{\begin{tabular}{|c|c|cc|ccc|} \toprule
estimator & coverage & neg.~bias & pos.~bias & lo.~quart. & median & up.~quart. \\ \colrule 
$\bar{y}^{*}(0.6092)$          & $100\%$ & $0\%$ & $0\%$ & $0.5732$ & $0.5829$ & $0.5928$\rule[-1.8mm]{0mm}{0mm} \\ 
$\bar{y}^{(-0.6092)}$             & $100\%$ & $0\%$ & $0\%$ & $0.5738$ & $0.5835$ & $0.5935$\rule[-1.8mm]{0mm}{0mm} \\ 
$\bar{y}_{0.6092}$               & $100\%$ & $0\%$ & $0\%$ & $0.5732$ & $0.5829$ & $0.5928$\rule[-1.8mm]{0mm}{0mm} \\ 
$\bar{y}_{0.6092\,t}$            & $100\%$ & $0\%$ & $0\%$ & $0.5719$ & $0.5816$ & $0.5915$\rule[-1.8mm]{0mm}{0mm} \\ 
$\bar{y}_{RP}(0.8046)$         & $99.99\%$ & $0\%$ & $0.01\%$ & $0.5749$ & $0.5850$ & $0.5953$\rule[-1.8mm]{0mm}{0mm} \\ \botrule
\end{tabular}}
\begin{minipage}[c]{\textwidth}
\begin{minipage}[c]{.0015\textwidth}
\mbox{ }
\end{minipage}
\begin{minipage}[c]{.4\textwidth}
\tablefont
{\begin{tabular}{|c|cc|} \toprule
estimator & MSE & $\frac{\mse(\bar{y}^{}_{})}{\mse(\text{est})}$\rule[-2.2mm]{0mm}{0mm}\rule[3.8mm]{0mm}{0mm} \\ \colrule 
$\bar{y}^{*}(0.6092)$          & $0.000212$ & $598.93\%$\rule[-1.8mm]{0mm}{0mm}\\ 
$\bar{y}^{(-0.6092)}$            & $0.000214$ & $593.91\%$\rule[-1.8mm]{0mm}{0mm}\\ 
$\bar{y}_{0.6092}$               & $0.000212$ & $598.92\%$\rule[-1.8mm]{0mm}{0mm}\\
$\bar{y}_{0.6092\,t}$           & $0.000215$ & $590.51\%$\rule[-1.8mm]{0mm}{0mm}\\ 
$\bar{y}_{RP}(0.8046)$        & $0.000228$ & $555.09\%$\rule[-1.8mm]{0mm}{0mm}\\ \botrule
\end{tabular}}
\end{minipage}\hfill
\begin{minipage}[c]{.52\textwidth}
\includegraphics[scale=0.45]{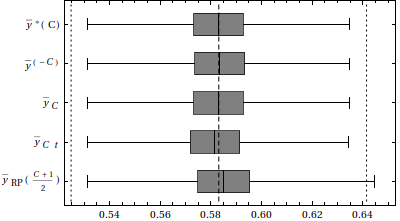}
\end{minipage}
\end{minipage}
\label{tab:est2}
\end{table}

\end{document}